\documentclass[12pt]{article}
\usepackage{amsmath,amssymb,amsthm,amsfonts}
\usepackage[T1]{fontenc}
\usepackage{color}

\newtheorem{thm}{Theorem}

\newtheorem{lem}[thm]{Lemma}

\theoremstyle{definition}
\newtheorem*{examp}{Example}

\theoremstyle{remark}
\newtheorem{rem}[thm]{Remark}

\newcommand{\thmref}[1]{Theorem~\ref{#1}}

\newcommand{\lemref}[1]{Lemma~\ref{#1}}

\newcommand{\abs}[1]{\left\lvert#1\right\rvert}

\newcommand{\cc}{\boldsymbol{\cdot}}
\newcommand{\Rn}{\mathbb{R}^n}
\newcommand{\RR}{\mathbb{R}}
\newcommand{\RRR}{{\mathbb{R}}^{n\times n}_{sym}}

\newcommand{\WWW}{H^1_0(\Omega;\RR^n)}

\newcommand{\ovm}{\overline{m}}
\newcommand{\unm}{\underline{m}}

\newcommand{\J}{\mathcal{J}}
\newcommand{\I}{\mathcal{I}}

\DeclareMathOperator{\Ker}{Ker}

\DeclareMathOperator{\const}{const}

\begin{document}


\par
\noindent \textbf{Nonconvex minimization
related to quadratic double-well energy -- approximation by convex problems}

\par
\vspace{1pc}
\noindent \framebox[4cm][c]{Zdzis{\l}aw Naniewicz}
\vspace{1pc}
\par
\noindent {\footnotesize
Cardinal Stefan Wyszy\'nski University\\
Faculty of Mathematics and Science\\
Dewajtis 5, 01-815 Warsaw, Poland \\
e-mail: naniewicz@uksw.edu.pl}

\par
\vspace{1pc}
\noindent Piotr Pucha{\l}a
\vspace{1pc}
\par
\noindent {\footnotesize
Cz\c{e}stochowa University of Technology\\
Faculty of
Mechanical Engineering and Computer Science\\
Institute of Mathematics\\
Armii Krajowej 21, 42-200 Cz\c{e}stochowa, Poland\\
e-mail: p.st.puchala@gmail.com, piotr.puchala@im.pcz.pl}

\vspace{1pc}
\begin{quotation}\noindent {\footnotesize \textbf{Abstract.}
A double-well energy expressed as a minimum of two
quadratic functions, called phase energies, is studied with taking into
account the minimization of the
corresponding integral functional. Such integral, as being not sequentially
weakly lower
semicontinuous, does not admit classical minimizers. To derive the relaxation
formula for the infimum, the minimizing sequence consisting of solutions
of convex problems appropriately approximating the original nonconvex one is
constructed.
The weak limit of this sequence together with the weak limit of the sequence
of solutions of the
corresponding dual problems and the weak limits of the characteristic
functions related to the phase energies are involved in the relaxation
formula.}
\end{quotation}
\begin{quotation}
\noindent {\footnotesize \textbf{Key Words.} Nonconvex minimization problem, Minimum of
convex functions, Duality, Pa\-ra\-me\-tri\-zed Young measures, Phase transitions.}
\end{quotation}
\begin{quotation}
\noindent {\footnotesize \textbf{AMS Classification.} 49J40, 49M20.}
\end{quotation}
\par

\renewcommand{\theequation}{\thesection.\arabic{equation}}
\section{Introduction}
\setcounter{section}{1}%
\setcounter{equation}{0}%

The minimization problem
of the form
\[
\inf\limits_{v\in \WWW}
\int_{\Omega}\min\Bigl\{\tfrac{1}{2}a
\abs{\varepsilon (v)+C}^2,\,
\tfrac{1}{2}b\abs{\varepsilon (v)+D}^2\Bigr\}\,dx \leqno{(P)}
\]
is considered, where $\varepsilon(v)$ is the symmetrized gradient of
$v\in\WWW$  and $\abs{\cdot }^2 $ stands for the square of
the scalar product ,,$\cdot\,$'' in
$L^2(\Omega;\RR^{n\times n}_{sym})$.\\

 Since the integrand involved is not quasiconvex, we cannot expect
the existence of classical minimizers of this functional.
\par
Basically, there are two ways to proceed if there are no minimizers as a
consequence of the lack of quasiconvexity (\cite{MORREY}). The first one is to
``quasiconvexify" the original functional and relate the gathered
information with the functional itself (cf. \cite{MORREY}, \cite{BALL},
\cite{BALL-MURAT}, \cite{ACERBI-FUSCO}, \cite{DACOROGNA}, \cite{KOHN},
\cite{KOHN-STRANG}, \cite{TARTAR-PREPRINT}, \cite{MURAT}, \cite{TARTAR-PROC},
\cite{FONSECA}, \cite{FONSECA-MULLER}, \cite{FONSECA-RYBKA}, \cite{BUTTAZZO},
\cite{DALMASO}, \cite{AMBROSIO}, \cite{BOUCHITTE-BRAIDES-BUTTAZZO},
\cite{ALLAIRE-FRANCFORT}, \cite{ALLAIRE-LODS} and the references quoted
there).  Another possibility is to enlarge the space of
admissible functions from Sobolev spaces to the space of parameterized
measures (called Young measures \cite{Young-37}) instead
of replacing the objective by its suitable envelopes. In this approach the
Young measures can be regarded
as means of summarizing the spatial oscillatory properties of minimizing
sequences, thus conserving some of that information.
With this respect we refer the reader to \cite{YOUNG},
\cite{KINDER-PEDREGAL}, \cite{CHIPOT-KINDER}, \cite{BALL-JAMES},
\cite{JAMES-KINDER},
\cite{BALL-MURAT}, \cite{MURAT}, \cite{ERICKSEN-PHASE}, \cite{PEDREGAL},
\cite{TARTAR-Y-H}  and the references therein.
\par
However, some important information is lost when seeking minima of
lower semicontinuous regularizations (quasiconvex envelopes). Minimizers
themselves are not sufficient
to characterize properly oscillatory phenomena of the
problem (microstructural features describing fine mixtures
of the phases in the phase transition problems, for instance).
From the application point of view, the detailed structure of minimizing
sequences appears to be as important as the minimizers
themselves.
Moreover, in the vectorial case it is almost impossible to
compute for a given objective its quasiconvexification. It is also very
difficult to compute explicitly the parametrized measures associated to a
minimizing sequence characterizing the infimum of the problem under
consideration.
\par
In this approach a new method to derive the formula for the infimum of $(P)$
is presented.
It preserves all the important information concerning the
oscillatory phenomena and is much easier to be obtained in practice.
The idea is to approximate $(P)$ by convex problems as proposed in
\cite{NAN-MINCONV}.
\begin{thm}[\cite{NAN-MINCONV}]\label{thm-1}
Suppose that $f_i:\Omega\times\mathbb{R}^m\times\mathbb{R}^{nm}\to
\mathbb{R}$, $i=1,2$, are quasiconvex, satisfy the Caratheodory and
growth conditions
\begin{description}
\item[(1)]  $\forall\,(s,\xi)\in{\mathbb R}^m\times
{\mathbb  R}^{nm}$,\quad $\Omega\ni x\mapsto f_i(x,s,\xi)$ is measurable,
\item[(2)] for a.e. $x\in \Omega,\quad {\mathbb R}^m\times
{\mathbb  R}^{nm}\ni (s,\xi)\mapsto f_i(x,s,\xi)$ is continuous,
\item[(3)] $a(x)+c(\abs{s}^2+\abs{\xi}^2)\leq f_i(x,s,\xi)\leq
A(x)+C(\abs{s}^2+\abs{\xi}^2)$,
\end{description}
where $a(\cdot)$, $A(\cdot)$ are non-negative summable functions in $\Omega$,
$c$ and $C$ are positive constants. Set
\[
\inf\biggl\{\int\limits_{\Omega}\min\bigl\{
f_1\bigl(x,u(x),\nabla u(x)\bigr),
f_2\bigl(x,u(x),\nabla u(x)\bigr)
\bigr\}\,dx\colon\; u\in \WWW\biggr\}:
=\alpha.                    \leqno{(Q)}
\]
Then there exist sequences
$u_k \in H^1_0(\Omega;{\mathbb R}^m)$, $\chi^{k}_1\colon\Omega\to\{0,1\}$ and
$\chi^{k}_2\colon\Omega\to\{0,1\}$ with $\chi^{k}_1+\chi^{k}_2\equiv 1$,
such that
\begin{description}
\item[$(i)$] $\{u_k\}$ is a minimizing sequence for $(Q)$,
\item[$(ii)$] $u_k\to u$ weakly in $H^1_0(\Omega;{\mathbb R}^m)$ as $k\to\infty$,
\item[$(iii)$] ${\chi}^{k}_1\to\chi_1$, ${\chi}^{k}_1\to\chi_2$
weak$^{\,\star}$ in $L^{\infty}(\Omega)$ as $k\to\infty$,
$\chi_1\colon\Omega\to[0,1]$, $\chi_2\colon\Omega\to[0,1]$ with
$\chi_1+\chi_2\equiv 1$,
\item[$(iv)$] $\lim\limits_{k\to\infty}\int\limits_{\Omega}\bigl[\chi^{k}_1 f_1(u_k)+
\chi^{k}_2f_2(u_k)\bigr]d\Omega=\alpha$,
\item[$(v)$] $\int\limits_{\Omega}\bigl[\chi^{k}_1 f_1(u_k)+
{\chi}^{k}_2f_2(u_k)\bigr]d\Omega
\leq\int\limits_{\Omega}\bigl[\chi^{k}_1f_1(w)+
\chi^{k}_2f_2(w)\bigr]d\Omega,\quad
\forall\,w\in H^1_0(\Omega;{\mathbb R}^m)$.
\end{description}
\end{thm}
The minimizing sequence for the problem $(P)$, established according to
\thmref{thm-1},
together with the sequence of solutions of the corresponding dual problems
and the sequence of characteristic functions related to the phases
$\tfrac{1}{2}a\abs{\varepsilon (\cdot)+C}^2$ and $\tfrac{1}{2}b\abs{
\varepsilon (\cdot)+D}^2$ generate limits (in appropriate weak topologies)
which are involved in the infimum formula.
It can be shown that the very special structure of the
minimizing sequence to be constructed allows the infimum to
be fully expressed by the parametrized
Young measures associated to this sequence. Some relations between the related
parametrized Young measures
and the weak limits of the characteristic functions are established.

\section{Statement of the problem and its approximation}

\renewcommand{\theequation}{\thesection.\arabic{equation}}
\setcounter{section}{2}%
\setcounter{equation}{0}%

Let $\Omega\subset\mathbb{R}^n$ be a bounded domain in $\mathbb{R}^n$ with
sufficiently smooth boundary $\partial \Omega$. Set
\[
\J(u)=\int\limits_{\Omega}\min\bigl\{\tfrac{1}{2}a
\abs{\varepsilon (u)+C}^2,
\tfrac{1}{2}b\abs{\varepsilon (u)+D}^2\bigr\}dx,\quad u\in\WWW.
\]
The problem to be considered here is
\[
\inf\Bigl\{\J(u)\colon \;u\in \WWW\Bigr\}:
=\alpha,                                         \leqno{(P)}
\]
where  $u:\Omega\subset \RR^n\to \RR^n$ is a competing vector-valued function
from the Sobolev space $\WWW$,
 $\varepsilon (u)\in
L^2(\Omega;\mathbb{R}^{n\times n}_{sym})$ is the symmetrized gradient of
$u\in\WWW$,
$C,D\in L^\infty(\Omega;\RR^{n\times n}_{sym})$ and where $a,b\in
L^\infty(\Omega)$ are such that $a(x),b(x)\geq\delta>0$ a.e. in $\Omega$ for a
positive constant $\delta$.
\par
\thmref{thm-1} ensures the existence of sequences $u^{k}\in \WWW$,
$\chi_a^{k}\subset\{0,1\}$
and $\chi_b^{k}\subset\{0,1\}$, $\chi_a^{k}+\chi_b^{k}\equiv 1$, with the
properties that
\begin{description}
\item[$(a)$] $\{u_k\}$ is a minimizing sequence for $(P)$,
\item[$(b)$] $u^{k}\to u$ weakly in $\WWW$ as $k\to\infty$,
\item[$(c)$] ${\chi}^{k}_a\to\chi_a,\quad {\chi}^{k}_b\to\chi_b$
$\text{weak}^
{\,\star}$ in $ L^{\infty}(\Omega)$ as $k\to\infty$, where
$\chi_a\colon\Omega\to[0,1]$,
$\chi_b\colon\Omega\to[0,1]$ with $\chi_a+\chi_b\equiv 1$,
\item[$(d)$] $\int\limits_\Omega
[\frac{1}{2}\chi_a^{k}a\abs{\varepsilon (u^{k})+C}^2
+\frac{1}{2}\chi_b^{k}b\abs{\varepsilon(u^{k})+D}^2]\, dx
:=\alpha^k\to \alpha$ as $k\to\infty$,
\item[$(e)$]
$
\int\limits_\Omega
[\frac{1}{2}\chi_a^{k}a\abs{\varepsilon (u^{k})+C}^2
+\frac{1}{2}\chi_b^{k}b\abs{\varepsilon (u^{k})+D}^2]\,dx\leq
$\\ 
\begin{flushright} $\leq\int\limits_\Omega
[\frac{1}{2}\chi_a^{k}a\abs{\varepsilon (w)+C}^2
+\frac{1}{2}\chi_b^{k}b\abs{\varepsilon (w)+D}^2]\,dx$,
$\quad\forall\;w\in\WWW$.
\end{flushright}
\end{description}

 Let us now introduce the function which describes the behaviour of the
minimizing sequence $\{ u_k\}$:\\
\begin{equation}\label{eq-1.1} 
\psi^k=\chi^k_b-\chi^k_a 
\end{equation}
with the property
\begin{equation}\label{eq-1.2}
(\psi^k)^2=1.
\end{equation}

It will be convenient to introduce the notations:
\begin{equation*}
m^k:=\chi_a^{k}a+\chi_b^{k}b,\quad
\ovm:=\frac{a+b}{2},\quad\unm:=\frac{b-a}{2}
\end{equation*}
Notice that $m^k$ has a decomposition
\begin{equation}
m^k=\overline{m}+\psi^k\underline{m}.
\end{equation}

From (e) of \thmref{thm-1} it follows that $u^k$ is a solution of the convex
optimization problem
\[
\inf\Bigl\{\J^k(v)\colon v\in \WWW\Bigr\}:=\alpha^k,   \leqno{(P^k)}
\]
where
\begin{equation*}
\J^k(v)=\int\limits_\Omega
\Bigl[\frac{1}{2}\chi_a^{k}a\abs{\varepsilon (v)+C}^2
+\frac{1}{2}\chi_b^{k}b\abs{\varepsilon (v)+D}^2\Bigr]\,dx,\quad v\in \WWW,
\end{equation*}
i.e.
\[
\J^k(u^k)=\alpha^k.
\]
Easy calculations using the properties of the scalar product
show that
\[
\chi_a^{k}aC+\chi_b^{k}bC=\frac{aC+bD}{2}+\psi^k\frac{bD-aC}{2},
\]
\[
\tfrac{1}{2}\chi^k_aa\abs{C}^2+
\tfrac{1}{2}\chi^k_bb\abs{D}^2=
\tfrac{1}{2}\biggl(\frac{a\abs{C}^2+b\abs{D}^2}{2}+\psi^k\frac{b\abs{D}^2
-a\abs{C}^2}{2}\biggr)
\]
so that
$\J^k(\cdot)$ admits the representation
\begin{equation} \tag{R}
\left.
\begin{array}{r}
\displaystyle \J^k(v)=
\int\limits_{\Omega}\Bigl[\tfrac{1}{2}(\chi^k_aa+
\chi^k_bb)\abs{\varepsilon (v)}^2
+(\chi^k_a aC+\chi^k_b bD)\cc\varepsilon (v)
\\
+\tfrac{1}{2}\chi^k_aa\abs{C}^2+
\tfrac{1}{2}\chi^k_bb\abs{D}^2\Bigr]dx
\\
\displaystyle=\int\limits_{\Omega}\left[\frac{1}{2}m^k\abs{\varepsilon (v)}^2+
(\mathcal{A}^++\psi^k \mathcal{A}^-)\cc\varepsilon (v)
+\frac{1}{2}\mathcal{B}^k\right]dx,
\end{array}
\right\}
\end{equation}

where the following notations have been introduced
\begin{align}\label{eq-0}
\begin{split}
\mathcal{A}^+&:=\frac{aC+bD}{2}
\\
\mathcal{A}^-&:=\frac{bD-aC}{2}
\\
\mathcal{B}^k&:=\frac{a\abs{C}^2+b\abs{D}^2}{2}+\psi^k\frac{b\abs{D}^2
-a\abs{C}^2}{2}.
\end{split}
\end{align}
\par
Now we associate to the minimization problem $(P^k)$ the dual problem
$(P^k)^\star$ following the idea of Fenchel \cite{FENCHEL} (cf. also
\cite{EKELAND-TEMAM}, \cite{AUBIN}).  Let $\phi(\cdot)$ be given by
\begin{equation}\label{eq-1.4}
\phi(\xi)=\frac{m}{2}\abs{\xi}^2+
E\cc\xi,\quad\xi\in\RRR,
\end{equation}
where $m>0$ is a constant,
then its conjugate $\phi^c(\cdot)$ reads
\begin{equation}\label{eq-1.5}
\phi^c(p)=\sup_{\xi\in\RRR}\bigl(p\cc\xi-\phi(\xi)\bigr)
=\frac{1}{2m}\abs{p-E}^2, \quad p\in\RRR.
\end{equation}
Next define a linear continuous operator $L:\WWW\to
L^2(\Omega;\mathbb{R}^{n\times n}_{sym})$ as
\[
Lv=\varepsilon(v),\quad v\in\WWW,
\]
then its transpose $L^\star:L^2(\Omega;\mathbb{R}^{n\times n}_{sym})\to
H^{-1}(\Omega;\mathbb{R}^n)$ can be expressed as
\[
\bigl\langle L^\star p,v\bigr\rangle_{\WWW}=\int\limits_\Omega p\cc\varepsilon
(v)\,dx,\quad \forall p\in L^2(\Omega;\mathbb{R}^{n\times n}_{sym}),\;\forall
v\in\WWW.
\]
Define
\begin{equation}\label{eq-1.66}
\I^k(q):=\int\limits_{\Omega}\Bigl[\frac{1}{2m^k}
\abs{q-(\mathcal{A}^++\psi^k\mathcal{A}^-)}^2-\frac{1}{2}\mathcal{B}^k\Bigr]
\,dx,\quad q\in
L^2(\Omega;\mathbb{R}^{n\times n}_{sym})
\end{equation}
and denote by $\Ker L^\star$ the kernel of $L^\star$, i.e.
\[
\Ker L^\star=\Bigl\{p\in L^2(\Omega;\mathbb{R}^{n\times n}_{sym})\colon
\int\limits_\Omega p\cc\varepsilon(v)\,dx=0\Bigl\}.
\]
\par
Now we are in a position to formulate the dual problem $(P^k)^\star$ which
can be stated as follows
\[
\inf\Bigl\{\int\limits_{\Omega} \I^k(q)\,dx\colon q\in\Ker
L^\star\Bigr\}:=\beta^k.                         \leqno{(P^k)^\star}
\]
According to the Fenchel theorem (cf. Theorem 3.2, p. 38, \cite{AUBIN})
we get
\begin{multline}\label{eq-1.6}
\J^k(v)=\int\limits_{\Omega}\Bigl[\frac{1}{2}m^k\abs{\varepsilon (v)}^2+
(\mathcal{A}^++\psi^k \mathcal{A}^-)\cc\varepsilon (v)
+\frac{1}{2}\mathcal{B}^k\Bigr]\,dx\geq
\\
\geq\,\int\limits_{\Omega}\Bigl[\frac{1}{2}m^k\abs{\varepsilon
(u^k)}^2+
(\mathcal{A}^++\psi^k \mathcal{A}^-)\cc\varepsilon (u^k)
+\frac{1}{2}\mathcal{B}^k\Bigr]\,dx=
\\
=\,\J^k(u^k)=\alpha^k=-\beta^k=-\I^k(p^k)=
\\
=-\int\limits_{\Omega}\Bigl[
\frac{1}{2m^k}
\abs{p^k-(\mathcal{A}^++\psi^k\mathcal{A}^-)}^2-\frac{1}{2}\mathcal{B}^k\Bigr]\,dx \geq
\\
\geq\, -\int\limits_{\Omega}\Bigl[\frac{1}{2m^k}
\abs{q-(\mathcal{A}^++\psi^k\mathcal{A}^-)}^2-\frac{1}{2}\mathcal{B}^k\Bigr]\,dx=-\I^k(q),
\\
\forall\,v\in \WWW,\;\forall\,q\in \Ker L^\star,
\end{multline}
where
\begin{equation}\label{1-10}
p^k=m^k\varepsilon (u^k)+\mathcal{A}^++\psi^k \mathcal{A}^-\in\Ker L^\star
\end{equation}
is a solution of the dual problem  $(P^k)^\star$. Since $p^k\in\Ker L^\star$,
\begin{equation}\label{1-1}
\int\limits_{\Omega}p^k\cc\varepsilon(v)\,dx=0,\quad\forall\,v\in\WWW,
\end{equation}
so, in particular,
\begin{equation}\label{1-2}
\int\limits_{\Omega}p^k\cc\varepsilon(u^k)\,dx=
\int\limits_{\Omega}\Bigl[m^k\abs{\varepsilon (u^k)}^2
+(\mathcal{A}^++\psi^k \mathcal{A}^-)\cc\varepsilon (u^k)\Bigr]dx=0.
\end{equation}

Taking into account that $\alpha^k=\I^k(p^k)$ and equation (\ref{1-2}) we get

the following representations 
\begin{equation}\label{eq-11}
\alpha^k=\tfrac{1}{2}\int\limits_{\Omega}\bigl[-m^k\abs{\varepsilon (u^k)}^2
+\mathcal{B}^k\bigr]dx
=\tfrac{1}{2}\int\limits_{\Omega}\bigl[(\mathcal{A}^++\psi^k \mathcal{A}^-)\cc\varepsilon (u^k)
+\mathcal{B}^k\bigr]dx.
\end{equation}
Analogously, from {{\eqref{1-10}}} we have 
$\varepsilon(u^k)=\frac{1}{m^k}(p^k-\mathcal{A}^+-\psi^k
\mathcal{A}^-)$, so from (\ref{1-2}) and \eqref{eq-11}
it follows
\[
\int\limits_{\Omega}\left[\frac{1}{m^k}\abs{p^k}^2
-\frac{1}{m^k}(\mathcal{A}^++\psi^k \mathcal{A}^-)\cc p^k\right]dx=0.
\]
Therefore
\begin{align}
\alpha^k&=-\tfrac{1}{2}\int\limits_{\Omega}\left[-\frac{1}{m^k}
\abs{p^k}^2+
\frac{1}{m^k}\abs{\mathcal{A}^++\psi^k \mathcal{A}^-}^2
-\mathcal{B}^k\right]dx=\notag
\\
&=\tfrac{1}{2}\int\limits_{\Omega}\left[\frac{1}{m^k}
(\mathcal{A}^++\psi^k \mathcal{A}^-)\cc p^k-
\frac{1}{m^k}\abs{\mathcal{A}^++\psi^k \mathcal{A}^-}^2 \label{eq-1}
+\mathcal{B}^k\right]dx
\end{align}
and we are led to the equality
\begin{equation}\label{eq-111}
\int\limits_{\Omega}\left[m^k\abs{\varepsilon (u^k)}^2+
\frac{1}{m^k}
\abs{p^k}^2\right]dx
=\int\limits_{\Omega}\left[
\frac{1}{m^k}\abs{\mathcal{A}^++\psi^k \mathcal{A}^-}^2\right]dx.
\end{equation}
\par
Now observe that the right hand side of this equality can be reorganized in
such a way that its limit as $k\to\infty$ is easy to be calculated. Indeed,
if we set
\[
\ovm^\natural:=\frac{1}{2}\bigl(\frac{1}{a}+\frac{1}{b}\bigr),\quad
\unm^\natural:=\frac{1}{2}\bigl(\frac{1}{b}-\frac{1}{a}\bigr),
\]
and recall that $\chi_{a}^{k}+\chi_{b}^{k}=1$ and
$(\psi^k)^2=1$ we obtain
\[
m^k =\frac{\ovm^\natural -\psi^k\unm^\natural} 
{(\ovm^\natural)^2 - (\unm^\natural)^2}
\]
so that
\[
\frac{1}{m^k}=\ovm^\natural+\unm^\natural\psi^k
\]
and from \eqref{eq-111}, by making use of \eqref{eq-1.2}, we obtain
\begin{equation}
\begin{array}{l}
\displaystyle \lim\limits_{k\to\infty}\int\limits_{\Omega}\left[m^k\abs{\varepsilon
(u^k)}^2+ \frac{1}{m^k}
\abs{p^k}^2\right]dx
\\
\displaystyle=\lim\limits_{k\to\infty}\int\limits_{\Omega}\left[
(\ovm^\natural+\unm^\natural\psi^k)\abs{\mathcal{A}^++\psi^k \mathcal{A}^-}^2\right]dx
\\
\displaystyle=\int\limits_{\Omega}\left[
(\ovm^\natural+\unm^\natural\,\psi)\abs{\mathcal{A}^+}^2
+2(\ovm^\natural\,\psi+\unm^\natural) \mathcal{A}^+\cc \mathcal{A}^-
+(\ovm^\natural+\unm^\natural\,\psi)\abs{\mathcal{A}^-}^2\right]dx
\\
\displaystyle=\tfrac{1}{2}\int\limits_{\Omega}\left[(\ovm^\natural+\unm^\natural\,\psi)
\bigl(a^2\abs{C}^2+b^2\abs{D}^2\bigr)+(\ovm^\natural\,\psi+\unm^\natural)
\bigl(b^2\abs{D}^2-a^2\abs{C}^2\bigr)
\right]dx
\\
\displaystyle=\tfrac{1}{2}\int\limits_{\Omega}\Bigl(a\abs{C}^2+b\abs{D}^2
\Bigr)dx+\tfrac{1}{2}\int\limits_{\Omega}\psi\Bigl(b\abs{D}^2-a\abs{C}^2\Bigr)
dx=\int\limits_{\Omega}\mathcal{B}\,dx,
\end{array}
\end{equation}
where
\[
\mathcal{B}=\frac{a\abs{C}^2+b\abs{D}^2}{2}+\psi\frac{b\abs{D}^2-a\abs{C}^2}{2}.
\]
\par
Let us introduce the set
\[
\Omega_0=\{x\in\Omega\colon a(x)=b(x)\}.
\]
Using \eqref{1-10} we obtain immediately that in $\Omega\setminus\Omega_0$,
\[
\psi^k\varepsilon(u^k)= \frac{1}{\unm}p^k-\frac{\ovm}{\unm}\varepsilon(u^k)-
\frac{1}{\unm}\mathcal{A}^+-\frac{1}{\unm}\psi^k \mathcal{A}^-.
\]
Thus from \eqref{eq-11} it follows
\begin{align}\label{eq-I}
\alpha^k=&\,\tfrac{1}{2}\int\limits_{\Omega\setminus\Omega_0}\bigl[
\mathcal{A}^+\cc\varepsilon (u^k)+
\mathcal{A}^-\cc\psi^k\varepsilon(u^k)+\mathcal{B}^k\bigr]dx \notag
\\
&+\tfrac{1}{2}\int\limits_{\Omega_0}\bigl[
\mathcal{A}^+\cc\varepsilon (u^k)+
\mathcal{A}^-\cc\psi^k\varepsilon(u^k)+\mathcal{B}^k\bigr]dx \notag
\\
=&\,\tfrac{1}{2}\int\limits_{\Omega\setminus\Omega_0}\Bigl[(\mathcal{A}^+-
\tfrac{\ovm}{\unm}\mathcal{A}^-)\cc\varepsilon(u^k)
+\tfrac{1}{\unm}\mathcal{A}^-\cc p^k-\tfrac{1}{\unm}\mathcal{A}^-
\cc \mathcal{A}^+  \notag
\\
&-\tfrac{1}{\unm}\psi^k\,\abs{
\mathcal{A}^-}^2+\mathcal{B}^k\Bigr]dx \notag
\\
&+\tfrac{1}{2}\int\limits_{\Omega_0}\bigl[
\mathcal{A}^+\cc\varepsilon (u^k)+
p^k\cc\varepsilon(u^k)-a\abs{\varepsilon(u^k)}^2-\mathcal{A}^+\cc\varepsilon(u^k)
+\mathcal{B}^k\bigr]dx  \notag
\\
=&\,\tfrac{1}{2}\int\limits_{\Omega\setminus\Omega_0}\Bigl[
\frac{ab(C-D)}{b-a}
\cc\varepsilon(u^k)
+\frac{bD-aC}{b-a}\cc p^k+
\frac{ab(\abs{C}^2-\abs{D}^2)}{2(b-a)}\notag
\\
&-\psi^k\,\frac{ab\abs{C-D}^2}{2(b-a)}
\Bigr]dx \notag
\\
&+\tfrac{1}{2}\int\limits_{\Omega_0}\Bigl[
-a\abs{\varepsilon(u^k)}^2+\frac{a(\abs{C}^2+\abs{D}^2)}{2}+
\psi^k\frac{a(\abs{D}^2-\abs{C}^2)}{2}\Bigr]dx \notag
\\
&+\tfrac{1}{2}\int\limits_{\Omega_0}p^k\cc\varepsilon(u^k)dx.
\end{align}
Another representation of $\alpha^k$ can be derived from \eqref{1-10} if we
take into account
\eqref{eq-1},
$(\psi^k)^2=1$ and that in $\Omega_0$,
$\frac{1}{a}\bigl(\abs{\mathcal{A}^+}^2+
\abs{\mathcal{A}^-}^2+2\psi^k\mathcal{A}^+\cc
\mathcal{A}^-\bigr)=\mathcal{B}^k$:
\begin{align}\label{eq-II}
\alpha^k=&\,\tfrac{1}{2}\int\limits_{\Omega\setminus\Omega_0}
\Bigl[\frac{1}{m^k}(\mathcal{A}^++\psi^k\mathcal{A}^-)\cc p^k-
\frac{1}{m^k}\abs{\mathcal{A}^++\psi^k \mathcal{A}^-}^2+\mathcal{B}^k
\Bigr]dx  \notag
\\
&+\tfrac{1}{2}\int\limits_{\Omega_0}
\Bigl[\frac{1}{a}(\mathcal{A}^++\psi^k\mathcal{A}^-)\cc p^k-
\frac{1}{a}\abs{\mathcal{A}^++\psi^k \mathcal{A}^-}^2
+\mathcal{B}^k\Bigr]dx \notag
\\
=&\,\tfrac{1}{2}\int\limits_{\Omega\setminus\Omega_0}\bigl[
\mathcal{A}^+\cc\varepsilon (u^k)+
\mathcal{A}^-\cc\psi^k\varepsilon(u^k)+\mathcal{B}^k\bigr]dx \notag
\\
&+\tfrac{1}{2}\int\limits_{\Omega_0}
\Bigl[\frac{1}{a}\abs{p^k}^2-
\frac{1}{a}\bigl(\abs{\mathcal{A}^+}^2+\abs{\mathcal{A}^-}^2+2\psi^k\mathcal{A}^+\cc \mathcal{A}^-\bigr)
+\mathcal{B}^k\Bigr]dx \notag
\\
=&\,\tfrac{1}{2}\int\limits_{\Omega\setminus\Omega_0}\Bigl[
\frac{ab(C-D)}{b-a}
\cc\varepsilon(u^k)
+\frac{bD-aC}{b-a}\cc p^k+
\frac{ab(\abs{C}^2-\abs{D}^2)}{2(b-a)}\notag
\\
&-\psi^k\frac{ab\abs{C-D}^2}{2(b-a)}
\Bigr]dx \notag
\\
&+\tfrac{1}{2}\int\limits_{\Omega_0}
\frac{1}{a}\abs{p^k}^2dx -\tfrac{1}{2}\int\limits_{\Omega_0}
p^k\cc\varepsilon(u^k)dx.
\end{align}
\par

As it will be seen in the next section the formulas \eqref{eq-I} and
\eqref{eq-II} allow us to express $\alpha$ in terms of the weak limits: $u$,
$p$ and $\psi$.

\section{Weak convergence in $\boldsymbol{L^1(\Omega)}$}
\renewcommand{\theequation}{\thesection.\arabic{equation}}
\setcounter{section}{3}%
\setcounter{equation}{0}%

\begin{lem}\label{lem-1}
Let $\Omega\subset \Rn$ be a bounded domain in $\Rn$
with Lipschitz continuous boundary $\partial\Omega$.
Then
\begin{equation}\label{P-1}
p^k\cc\varepsilon(u^k)\to p\cc\varepsilon(u)\quad\mbox{weakly in}
\;L^1(\Omega).
\end{equation}
\end{lem}
\begin{proof}
Extend each function $u^k\in\WWW$ to all of $\Rn$
by setting
it equal to zero on $\Rn\setminus\Omega$. By regularity of the boundary
$\partial\Omega$ all of these extensions are elements of $H^1(\Rn;\Rn)$.
For an arbitrary $\varphi\in C^\infty(\Rn)$ we thus have $\varphi u^k\in
H^1(\Rn;\Rn)$ and $\varphi u^k\big|_\Omega\in\WWW$. We claim that
\begin{equation}\label{P-3}
\int\limits_{\Rn}\varphi \,p^k\cc\varepsilon(u^k)\,dx\to
\int\limits_{\Rn}\varphi \,p\cc\varepsilon(u)\,dx,
\end{equation}
for any $\varphi\in C^\infty(\Rn)$.
Indeed, since $\varepsilon(\varphi
u^k)=\varphi\varepsilon(u^k)+u^k\otimes\nabla \varphi$,
we get
\begin{eqnarray*}
\int\limits_{\Rn}\varphi \,p^k\cc\varepsilon(u^k)\,dx=
\int\limits_{\Rn}
p^k\cc\varepsilon(\varphi u^k)\,dx-\int\limits_{\Rn}
p^k\cc (u^k\otimes\nabla \varphi)\,dx=
\\
=-\int\limits_{\Rn}
p^k\cc (u^k\otimes\nabla \varphi)\,dx\to-\int\limits_{\Rn}
p\cc (u\otimes\nabla \varphi)\,dx=
\int\limits_{\Rn}\varphi
\,p\cc\varepsilon(u)\,dx,
\end{eqnarray*}
where we have used \eqref{1-1} and the strong convergence $u^k\to u$ in
$L^2(\Omega;\Rn)$ (valid due to the Rellich compactness theorem).
\par
Further, the sequence $\{p^k\cc\varepsilon(u^k)\}$ is uniformly bounded in
$L^1(\Omega)$
because
$\{p^k\}$ and $\{\varepsilon(u^k)\}$ so are in $L^2(\Omega)$. By
Chacon's
biting lemma \cite{PEDREGAL} it follows that there exist a subsequence of
$(p^k\cc\varepsilon(u^k))$, not relabeled, a
nonincreasing sequence of measurable sets $\Omega_n\subset\Omega$,
$\Omega_{n+1}\subset \Omega_n$,
$\abs{\Omega_n} \searrow 0$ and $f\in L^1(\Omega)$ such that
\begin{equation}\label{P-2}
p^k\cc\varepsilon(u^k)\to f\quad\mbox{weakly in}\;L^1(\Omega\setminus\Omega_n)
\end{equation}
for all $n$. It means that $\{p^k\cc\varepsilon(u^k)\}$ converges in the
biting sense to $f$ \cite{PEDREGAL}.
\par
Now we assert that the biting limit $f$ coincides with $p\cc\varepsilon(u)$,
i.e.  $f=p\cc\varepsilon(u)$ a.e. in $\Omega$. To show this
observe that from  the biting argument \eqref{P-2} and \eqref{P-3} it follows
that for any $\varphi\in C^\infty(\Rn)$ we get
\begin{equation}\label{P-4}
\int\limits_{\Omega\setminus\Omega_n}\varphi\,p\cc\varepsilon(u)\,dx=
\int\limits_{\Omega\setminus\Omega_n}\varphi\,f\,dx,
\end{equation}
for any $n$. Hence $p\cc\varepsilon(u)=f$ a.e. in
$\Omega\setminus\Omega_n$ for each $n$. Since $\abs{\Omega_n}\searrow 0$
as $n\to\infty$, the equality $p\cc\varepsilon(u)=f$ must hold a.e. in
$\Omega$. Thus the assertion follows.
\par
Recall that $p^k\cc\varepsilon(u^k)=m^k\abs{\varepsilon(u^k)}^2+
\bigl(\mathcal{A}^++\psi^k \mathcal{A}^-)\cc\varepsilon(u^k)$. Therefore one can deduce the
existence of a constant $C\geq 0$ such that
\[
p^k\cc\varepsilon(u^k)+C\geq 0\quad\mbox{a.e in }\Omega.
\]
Obviously $p^k\cc\varepsilon(u^k)+C$ converges in the biting sense to
$p\cc\varepsilon(u)+C$. According to (Lemma 6.9, p.109, \cite{PEDREGAL}) its
weak convergence in $L^1(\Omega)$ is then equivalent to
\begin{equation}\label{P-5}
\limsup\limits_{k\to\infty}\int\limits_\Omega
\bigl(p^k\cc\varepsilon(u^k)+C\bigr)\,dx\leq
\int\limits_\Omega \bigl(p\cc\varepsilon(u)+C\bigr)\,dx.
\end{equation}
Our task now is to establish the foregoing inequality. For this purpose notice
that \eqref{P-3} can be easily extend to the convergence
\begin{equation}\label{P-33}
\int\limits_{\Rn}\varphi \,p^k\cc\varepsilon(u^k)\,dx\to
\int\limits_{\Rn}\varphi \,p\cc\varepsilon(u)\,dx,
\end{equation}
which is valid for any $\varphi\in C_c(\Rn)$, where $C_c(\Rn)$ is the space of
continuous functions on $\Rn$ with compact support.  Thus
$\mu^k:=(p^k\cc\varepsilon(u^k)+C)dx$ and $\mu:=(p\cc\varepsilon(u)+C)dx$ can
be treated as positive Radon measures on $\Rn$ for which it holds
\[
\lim\limits_{k\to\infty}\int\limits_{\Rn}\varphi\,d\mu^k=
\int\limits_{\Rn}\varphi\,d\mu,\quad\forall\,\varphi\in C_c(\Rn).
\]
But (Theorem 1, p. 54, \cite{EVANS-GARIEPY}) asserts that this condition is
equivalent to the following one
\begin{equation}\label{P-6}
\lim\limits_{k\to\infty}\mu^k(B)=\mu(B)\quad \mbox{for each bounded Borel set
$B\subset\Rn$ with $\mu(\partial B)=0$}.
\end{equation}
Now we are in a position to show \eqref{P-5}. Fix $\epsilon>0$ and choose
$0<\delta<\epsilon$ with the property that $\omega\subset\Omega$ with
$\abs{\omega}<\delta$ implies
\[
\int\limits_\omega p\cc\varepsilon(u)\,dx<\epsilon.
\]
In the biting convergence take $n_0$ large enough to fulfill
$\abs{\Omega_{n_0}}<\tfrac{\delta}{2}$. By the measurability of $\Omega_{n_0}$
there exists an open $\widetilde{\Omega}_{n_0}\supset\Omega_{n_0}$ with
$\abs{\widetilde{\Omega}_{n_0}}<\delta$. Vitali's covering theorem ensures the
representation $\widetilde{\Omega}_{n_0}=\widetilde{\Omega}_{n_0}^\prime\cup
\widetilde{\Omega}_{n_0}^{\prime\prime}$ where
$\abs{\widetilde{\Omega}_{n_0}^{\prime\prime}}=0$ and
$\widetilde{\Omega}_{n_0}^\prime$ stands for the union of a countable
collection of disjoint closed balls in $\widetilde{\Omega}_{n_0}$. Therefore
$\abs{\partial \widetilde{\Omega}_{n_0}^\prime}=0$ and consequently
$\mu(\partial \widetilde{\Omega}_{n_0}^\prime)=0$. From this we have
\begin{align*}
\int\limits_\Omega \bigl(p^k\cc\varepsilon(u^k)+C\bigr)\,dx&=
\int\limits_{\Omega\setminus\Omega_{n_0}} \bigl(p^k\cc\varepsilon(u^k)+C\bigr)
\,dx+\int\limits_{\Omega_{n_0}}
\bigl(p^k\cc\varepsilon(u^k)+C\bigr)\,dx
\\
&\leq \int\limits_{\Omega\setminus\Omega_{n_0}}
\bigl(p^k\cc\varepsilon(u^k)+C\bigr) \,dx
+\int\limits_{\widetilde{\Omega}_{n_0}} \bigl(p^k\cc\varepsilon(u^k)+C\bigr)
\,dx
\\
&= \int\limits_{\Omega\setminus\Omega_{n_0}}
\bigl(p^k\cc\varepsilon(u^k)+C\bigr) \,dx+
\int\limits_{\widetilde{\Omega}_{n_0}^\prime}
\bigl(p^k\cc\varepsilon(u^k)+C\bigr) \,dx
\\
&=\int\limits_{\Omega\setminus\Omega_{n_0}}
\bigl(p^k\cc\varepsilon(u^k)+C\bigr) \,dx+
\mu^k(\widetilde{\Omega}_{n_0}^\prime)
\end{align*}
which thanks to \eqref{P-6} by passing to the limit as $k\to\infty$ yields
\begin{eqnarray*}
\limsup\limits_{k\to\infty}\int\limits_\Omega
\bigl(p^k\cc\varepsilon(u^k)+C\bigr)\,dx
\leq\int\limits_{\Omega\setminus\Omega_{n_0}} \bigl(p\cc\varepsilon(u)+C
\bigr)\,dx+
\mu(\widetilde{\Omega}_{n_0}^\prime)
\\
\leq\int\limits_{\Omega} \bigl(p\cc\varepsilon(u)+C\bigr)\,dx+
\int\limits_{\widetilde{\Omega}_{n_0}^\prime} \bigl(p\cc\varepsilon(u)+C
\bigr)\,dx
\\
\leq\int\limits_{\Omega} \bigl(p\cc\varepsilon(u)+C\bigr)\,dx+\epsilon(1+C),
\end{eqnarray*}
because $\abs{\widetilde{\Omega}_{n_0}^\prime}<\delta<\epsilon$.
Since $\epsilon>0$ was chosen arbitrarily, \eqref{P-5} follows. This completes
the proof of \lemref{lem-1}.
\end{proof}

\par
\section{Relaxed formulas for the infimum}
\renewcommand{\theequation}{\thesection.\arabic{equation}}
\setcounter{section}{4}%
\setcounter{equation}{0}%

The weak lower semicontinuity of convex functionals, the upper semicontinuity
of concave functionals and \lemref{lem-1} yield
\begin{equation}
\liminf\limits_{k\to\infty}\;\tfrac{1}{2}\int\limits_{\Omega_0}
\frac{1}{a}\abs{p^k}^2\,dx
\geq
\tfrac{1}{2}\int\limits_{\Omega_0}
\frac{1}{a}\abs{p}^2\,dx,
\end{equation}
\begin{equation}
\limsup\limits_{k\to\infty}\,\tfrac{1}{2}\int\limits_{\Omega_0}\Bigl[
-a\abs{\varepsilon(u^k)}^2+\mathcal{B}^k\Bigr]dx\leq
\tfrac{1}{2}\int\limits_{\Omega_0}\Bigl[
-a\abs{\varepsilon(u)}^2+\mathcal{B}\Bigr]dx,
\end{equation}
\begin{equation}
\lim\limits_{k\to\infty}\;\tfrac{1}{2}\int\limits_{\Omega_0}
p^k\cc\varepsilon(u^k)
dx=\tfrac{1}{2}\int\limits_{\Omega_0}p\cc\varepsilon(u)dx,
\end{equation}
where $\psi=\chi_b-\chi_a$ and
\[
\mathcal{B}=\frac{a(\abs{C}^2+\abs{D}^2)}{2}+\psi\frac{a(\abs{D}^2
-\abs{C}^2)}{2}\quad\mbox{in }\Omega_0.
\]
\par
Now we show that
\begin{eqnarray}
\lim\limits_{k\to\infty}\tfrac{1}{2}\int\limits_{\Omega\setminus\Omega_0}
\Bigl[
\frac{ab(C-D)}{b-a}
\cc\varepsilon(u^k)
+\frac{bD-aC}{b-a}\cc p^k+
\frac{ab(\abs{C}^2-\abs{D}^2)}{2(b-a)}\notag
\\
-\psi^k\frac{ab\abs{C-D}^2}{2(b-a)}
\Bigr]dx \notag
\\
=\tfrac{1}{2}\int\limits_{\Omega\setminus\Omega_0}\Bigl[
\frac{ab(C-D)}{b-a}
\cc\varepsilon(u)
+\frac{bD-aC}{b-a}\cc p+
\frac{ab(\abs{C}^2-\abs{D}^2)}{2(b-a)}\notag
\\
-\psi\frac{ab\abs{C-D}^2}{2(b-a)}
\Bigr]dx.   \label{eq-00}
\end{eqnarray}
This is not trivial because the functions $\frac{ab(C-D)}{b-a}$ and $\frac{bD-aC}{b-a}$ are not
assumed to belong to $L^2(\Omega\setminus\Omega_0;\mathbb{R}^{n\times n}_{sym})$. To overcome
this disadvantage let us recall (see \eqref{eq-I} or \eqref{eq-II})
that $\Omega\setminus\Omega_0$ is a set of a finite 
Lebesgue measure where
\begin{multline*}
\frac{ab(C-D)}{b-a}
\cc\varepsilon(u^k)
+\frac{bD-aC}{b-a}\cc p^k+
\frac{ab(\abs{C}^2-\abs{D}^2)}{2(b-a)}-\psi^k\,\frac{ab\abs{C-D}^2}{2(b-a)}
\\
=
\mathcal{A}^+\cc\varepsilon (u^k)+
\mathcal{A}^-\cc\psi^k\varepsilon(u^k)+\mathcal{B}^k.
\end{multline*}
Thus for any $\varepsilon >0$ there exist $\omega_\varepsilon\subset\Omega\setminus\Omega_0$
and $\delta>0$ such that $\abs{\omega_\varepsilon}<\varepsilon$ and for each $x\in
(\Omega\setminus\Omega_0)\setminus \omega_\varepsilon$ one has $\abs{a(x)-b(x)}\geq\delta$.
Hence
\begin{equation}
\frac{C-D}{b-a}\in L^\infty\bigl((\Omega\setminus\Omega_0)\setminus\omega_\varepsilon;
\mathbb{R}^{n\times n}_{sym}\bigr)\subset
L^2\bigl((\Omega\setminus\Omega_0)\setminus\omega_\varepsilon;
\mathbb{R}^{n\times n}_{sym}\bigr)
\end{equation}
and
\begin{multline*}
\Big\lvert\,\int\limits_{\omega_\varepsilon}
\Bigl[
\frac{ab(C-D)}{b-a}
\cc\varepsilon(u^k)
+\frac{bD-aC}{b-a}\cc p^k+
\frac{ab(\abs{C}^2-\abs{D}^2)}{2(b-a)}\notag
\\
-\psi^k\frac{ab\abs{C-D}^2}{2(b-a)}
\Bigr]dx\Big\rvert
\leq \const \abs{\omega_\varepsilon}^\frac{1}{2}\leq \const \varepsilon^\frac{1}{2}.
\end{multline*}
This allows the conclusion that
\begin{eqnarray*}
\lim\limits_{k\to\infty}
\int\limits_{(\Omega\setminus\Omega_0)\setminus\omega_\varepsilon}
\Bigl[
\frac{ab(C-D)}{b-a}
\cc\varepsilon(u^k)
+\frac{bD-aC}{b-a}\cc p^k+
\frac{ab(\abs{C}^2-\abs{D}^2)}{2(b-a)}
\\
-\psi^k\frac{ab\abs{C-D}^2}{2(b-a)}
\Bigr]dx
\\
=\int\limits_{(\Omega\setminus\Omega_0)\setminus\omega_\varepsilon}\Bigl[
\frac{ab(C-D)}{b-a}
\cc\varepsilon(u)
+\frac{bD-aC}{b-a}\cc p+
\frac{ab(\abs{C}^2-\abs{D}^2)}{2(b-a)}
\\
-\psi\frac{ab\abs{C-D}^2}{2(b-a)}
\Bigr]dx
\end{eqnarray*}
and due to the fact that $\varepsilon>0$ was chosen arbitrarily we easily arrive at
\eqref{eq-00}, as desired.
\par

Now, for $v\in \WWW$ and $q\in \Ker L^\star$ let us set
\begin{multline}
\mathcal{I}(v,q):=
\\
\int\limits_{\Omega\setminus\Omega_0}\Bigl[
\frac{ab(C-D)}{b-a}
\cc\varepsilon(v)
+\frac{bD-aC}{b-a}\cc q+
\frac{ab(\abs{C}^2-\abs{D}^2)}{2(b-a)}-\psi\,\frac{ab\abs{C-D}^2}{2(b-a)}
\Bigr]dx.
\end{multline}

Using the fact that for all $k\in \mathbb{N}$ we have
$p^k\in\Ker L^{\ast}$, \eqref{eq-00} and
passing to the limit as $k\to\infty$ in \eqref{eq-I} we
get 
\[
\alpha\leq\tfrac{1}{2}\mathcal{I}(u,p)+
\tfrac{1}{2}\int\limits_{\Omega_0}\Bigl[-a\abs{\varepsilon(u)}^2\Bigr]dx+
\tfrac{1}{2}\int\limits_{\Omega_0}\mathcal{B}\,dx+
\tfrac{1}{2}\int\limits_{\Omega_0}
p\cc\varepsilon(u)dx.
\]
Analogously, passing to the limit as $k\to\infty$ in  
\eqref{eq-II} we get
\[
\alpha\geq\tfrac{1}{2}\mathcal{I}(u,p)+
\tfrac{1}{2}\int\limits_{\Omega_0}
\Bigl[\frac{1}{a}\abs{p}^2-
\frac{1}{a}\bigl(\abs{\mathcal{A}^+}^2+\abs{\mathcal{A}^-}^2+2\psi
\mathcal{A}^+\cc \mathcal{A}^-\bigr) \Bigr]dx+
\tfrac{1}{2}\int\limits_{\Omega_0}\mathcal{B}\,dx.
\]
Putting the two above inequalities together we obtain the system of inequalities
\begin{multline}\label{eq-III}
\tfrac{1}{2}\int\limits_{\Omega_0}\Bigl[
-a\abs{\varepsilon(u)}^2\Bigr]dx
+\int\limits_{\Omega_0}p\cc\varepsilon(u)dx\geq
\\
\geq\alpha+\tfrac{1}{2}\int\limits_{\Omega_0}
p\cc\varepsilon(u)dx
-\mathcal{I}(u,p)-\int\limits_{\Omega_0}\mathcal{B}\,dx\geq
\\
\geq\tfrac{1}{2}\int\limits_{\Omega_0}
\Bigl[\frac{1}{a}\abs{p}^2-
\frac{1}{a}\bigl(\abs{\mathcal{A}^+}^2+\abs{\mathcal{A}^-}^2+2\psi
\mathcal{A}^+\cc \mathcal{A}^-\bigr) \Bigr]dx.
\end{multline}
But from the fact that $p=a\,\varepsilon(u)+\mathcal{A}^++\psi \mathcal{A}^-$
in $\Omega_0$ it follows
\begin{multline}
\tfrac{1}{2}\int\limits_{\Omega_0}
\Bigl[\frac{1}{a}\abs{p}^2-
\frac{1}{a}\bigl(\abs{\mathcal{A}^+}^2+\abs{\mathcal{A}^-}^2+2\psi \mathcal{A}^+\cc \mathcal{A}^-\bigr)
\Bigr]dx= \notag
\\
=\tfrac{1}{2}\int\limits_{\Omega_0}
\Bigl[a\abs{\varepsilon(u)}^2+\frac{1}{a}\abs{\mathcal{A}^++\psi\,\mathcal{A}^-}^2+
2(\mathcal{A}^++\psi\,\mathcal{A}^-)\cc\varepsilon(u)+\notag
\\
-\frac{1}{a}\bigl(\abs{\mathcal{A}^+}^2+\abs{\mathcal{A}^-}^2+2\psi \mathcal{A}^+\cc \mathcal{A}^-\bigr)
\Bigr]dx= \notag
\\
=\,\tfrac{1}{2}\int\limits_{\Omega_0}
\Bigl[-a\abs{\varepsilon(u)}^2+
2\bigl(a\abs{\varepsilon(u)}^2+(\mathcal{A}^++\psi\,\mathcal{A}^-)\cc\varepsilon(u)\bigl)+
\frac{1}{a}(\psi^2-1)\abs{\mathcal{A}^-}^2\Bigr]dx= \notag
\\
=\,\tfrac{1}{2}\int\limits_{\Omega_0}\Bigl[
-a\abs{\varepsilon(u)}^2\Bigr]dx
+\int\limits_{\Omega_0}p\cc\varepsilon(u)dx
-2\int\limits_{\Omega_0}\frac{1}{a}\,\chi_a\chi_b\,\abs{\mathcal{A}^-}^2dx.
\end{multline}
Here we used the fact that $\psi^2-1=-4\chi_a\chi_b$. Thus in view of
\eqref{eq-III} it follows
\begin{multline}
0\geq\alpha-\tfrac{3}{2}\int\limits_{\Omega_0}
p\cc\varepsilon(u)dx+\int\limits_{\Omega_0}
a\abs{\varepsilon(u)}^2dx
-\mathcal{I}(u,p)-\int\limits_{\Omega_0}\mathcal{B}\,dx\geq
\\
\geq -2\int\limits_{\Omega_0}\frac{1}{a}\,\chi_a\chi_b\,\abs{\mathcal{A}^-}^2
dx.
\end{multline}
Since
\[
\int\limits_{\Omega_0}p\cc\varepsilon(u)dx= \int\limits_{\Omega_0}\Bigl[
a\abs{\varepsilon(u)}^2+(\mathcal{A}^++\psi\,\mathcal{A}^-)\cc
\varepsilon(u)\Bigr]dx
\]
and in $\Omega_0$ $\mathcal{A}^- =a\frac{D-C}{2}$
we have
\begin{equation}
0\geq\alpha-\int\limits_{\Omega_0}\bigl[
(\mathcal{A}^++\psi\,\mathcal{A}^-)\cc\varepsilon(u)+\mathcal{B}\bigr]dx
-\mathcal{I}(u,p)-\tfrac{1}{2}\int\limits_{\Omega_0}p\cc\varepsilon(u)dx\\
\geq -\tfrac{1}{2}\int\limits_{\Omega_0}\chi_a\chi_b\,a\abs{C-D}^2dx.
\end{equation}
Thus we are allowed to conclude that
there exists a $\theta\in[0,1]$ such that
\begin{multline*}
\alpha=\,\int\limits_{\Omega_0}\bigl[
(\mathcal{A}^++\psi\,\mathcal{A}^-)\cc\varepsilon(u)+\mathcal{B}\bigr]dx
+\tfrac{1}{2}\int\limits_{\Omega_0}p\cc\varepsilon(u)dx
\\
+\int\limits_{\Omega\setminus\Omega_0}\Bigl[
\frac{ab(C-D)}{b-a}
\cc\varepsilon(u)
+\frac{bD-aC}{b-a}\cc p
+\frac{ab(\abs{C}^2-\abs{D}^2)}{2(b-a)}
\\
-\psi\,\frac{ab\abs{C-D}^2}{2(b-a)}
\Bigr]dx -\,\tfrac{\theta}{2}\int\limits_{\Omega_0}
\chi_a\chi_b \,a\abs{C-D}^2dx.
\end{multline*}
The obtained result can be summarized as follows.
\begin{thm}\label{thm-2}
Let $u\in\WWW$ and $p\in L^2(\Omega;\RRR)$ are the weak limits of $\{u^k\}$
and $\{p^k\}$ as defined by \eqref{eq-1.6}, respectively. Then there exists
$\theta\in[0,1]$ such that
\begin{align}\label{eq-IV}
\alpha=&\,\int\limits_{\Omega_0}\Bigl[
\Bigl(\frac{a(C+D)}{2}+\psi \frac{a(D-C)}{2}\Bigr)\cc\varepsilon(u)+
\frac{a(\abs{C}^2+\abs{D}^2)}{2}\notag
\\
&+\psi\frac{a(\abs{D}^2-\abs{C}^2)}{2}
\Bigr]dx+\tfrac{1}{2}\int\limits_{\Omega_0}p\cc\varepsilon(u)dx \notag
\\
&+\int\limits_{\Omega\setminus\Omega_0}\Bigl[
\frac{ab(C-D)}{b-a}
\cc\varepsilon(u)
+\frac{bD-aC}{b-a}\cc p
+\frac{ab(\abs{C}^2-\abs{D}^2)}{2(b-a)}\notag
\\
&-\psi\,\frac{ab\abs{C-D}^2}{2(b-a)}
\Bigr]\,dx-\,\tfrac{\theta}{2}\int\limits_{\Omega_0}
\chi_a\chi_b \,a\abs{C-D}^2dx.
\end{align}
\end{thm}

Using the fact that in $\Omega_0$ there holds the
equality
\[
\int\limits_{\Omega_0 }(\mathcal{A}^+ +\psi^k\mathcal{A}^-)
\varepsilon (u)dx=
\int\limits_{\Omega_0 }[-a\abs{\varepsilon (u)}^2]dx+
\int\limits_{\Omega_0 }p\cc\varepsilon (u)dx,
\]
we get from the above formula
\begin{align}
\alpha=&\,\int\limits_{\Omega_0}\Bigl[
-a\abs{\varepsilon(u)}^2+
\frac{a(\abs{C}^2+\abs{D}^2)}{2}+\psi\frac{a(\abs{D}^2-\abs{C}^2)}{2}\Bigr]dx
+\tfrac{3}{2}\int\limits_{\Omega_0}p\cc\varepsilon(u)dx+ \notag
\\
&+\,\int\limits_{\Omega\setminus\Omega_0}\Bigl[
\frac{ab(C-D)}{b-a}
\cc\varepsilon(u)
+\frac{bD-aC}{b-a}\cc p+
\frac{ab(\abs{C}^2-\abs{D}^2)}{2(b-a)}\notag
\\
&-\psi\frac{ab\abs{C-D}^2}{2(b-a)}
\Bigr]dx-\,\tfrac{\theta}{2}\int\limits_{\Omega_0}\chi_a\chi_b \,a\abs{C-D}^2dx
\label{eq-1.28}
\end{align}
and from \eqref{eq-1.28} and again the equality
$
\int\limits_{\Omega_0 }(\mathcal{A}^+ +\psi^k\mathcal{A}^-)
\varepsilon (u)dx=
\int\limits_{\Omega_0 }[-a\abs{\varepsilon (u)}^2]dx+\\
+\int\limits_{\Omega_0 }p\cc\varepsilon (u)dx,
$
\begin{align}
\alpha=
&\,\int\limits_{\Omega_0}
\frac{1}{a}\abs{p}^2\,dx -\tfrac{1}{2}\int\limits_{\Omega_0}
p\cc\varepsilon(u)\,dx+ \notag
\\
&+\,\int\limits_{\Omega\setminus\Omega_0}\Bigl[
\frac{ab(C-D)}{b-a}
\cc\varepsilon(u)
+\frac{bD-aC}{b-a}\cc p+
\frac{ab(\abs{C}^2-\abs{D}^2)}{2(b-a)}\notag
\\
&-\psi\frac{ab\abs{C-D}^2}{2(b-a)}
\Bigr]\,dx+ \tfrac{2-\theta}{2}\int\limits_{\Omega_0}\chi_a\chi_b \,
a\abs{C-D}^2dx.\label{eq-1.30}
\end{align}
By adding \eqref{eq-1.28} and \eqref{eq-1.30} we get also the formula:
\begin{align}
\alpha
=&\, {\tfrac{1}{2}}\int\limits_{\Omega_0}
\Bigl[\frac{1}{a}\abs{p}^2-a\abs{\varepsilon(u)}^2
+\frac{a(\abs{C}^2+\abs{D}^2)}{2}+\psi\frac{a(\abs{D}^2-\abs{C}^2)}{2}
\Bigr]dx+\tfrac{1}{2}\int\limits_{\Omega_0}p\cc\varepsilon(u)dx+ \notag
\\
&+\, \int\limits_{\Omega\setminus\Omega_0}\Bigl[
\frac{ab(C-D)}{b-a}
\cc\varepsilon(u)
+\frac{bD-aC}{b-a}\cc p+
\frac{ab(\abs{C}^2-\abs{D}^2)}{2(b-a)}\notag
\\
&-\psi\frac{ab\abs{C-D}^2}{2(b-a)}
\Bigr]dx+ {\tfrac{1-\theta}{2}}\int\limits_{\Omega_0}\chi_a\chi_b \,
a\abs{C-D}^2dx. \label{eq-1.31}
\end{align}
\newpage 
Before formulating next theorem it will be convenient to introduce
some notation. Denote by $\omega_0^+$ and $\omega_0^-$ such subsets of
$\Omega$, that $\psi^k\rightarrow 1$ weakly in $L^1(\omega_0^+)$ and 
$\psi^k\rightarrow -1$ weakly in $L^1(\omega_0^-)$. Let
$\omega_0\colon =\omega_0^+\cup\omega_0^-$.
\begin{thm}\label{thm-3}
Let $\nu=\{\nu_x\}_{x\in\Omega}$ be the  parametrized
Young measure associated to the mi\-ni\-mi\-za\-tion sequence
$\{u^k\}$. Then
\begin{align}\label{P-7}
\alpha =&\int\limits_{\Omega}\int\limits_{\Rn}h(x,\lambda)\,
d\nu_x(\lambda)\,dx=\notag
\\
=&\,\int\limits_{\Omega\setminus\Omega_0}\Bigl[
\frac{ab(C-D)}{b-a}
\cc\varepsilon(u)
+\frac{bD-aC}{b-a}\cc p+
\frac{ab(\abs{C}^2-\abs{D}^2)}{2(b-a)}\notag
\\
&-\psi\,\frac{ab\abs{C-D}^2}{2(b-a)}
\Bigr]dx \notag
\\
&+\int\limits_{\Omega_0}\Bigl[
-\int\limits_{\mathbb{R}^{{{n\times n}}}}a
{{|}}\lambda {{|}}^2\,d\nu_x(\lambda)+
\frac{a(\abs{C}^2+\abs{D}^2)}{2}+\psi\frac{a(\abs{D}^2-\abs{C}^2)}{2}\Bigr]dx\notag
\\
&+\tfrac{3}{2}\int\limits_{\Omega_0}p\cc\varepsilon(u)dx,
\end{align}
where $h(x,\lambda)=\min\bigl\{\tfrac{1}{2}a(x)
\abs{\lambda+C(x)}^2,\,
\tfrac{1}{2}b(x)\abs{\lambda+D(x)}^2\bigr\}$, $\lambda\in\RRR$, $x\in\Omega$.
Moreover, we have
\begin{equation}\label{P-71}
\nu_x=\delta_{\varepsilon(u(x))}\quad\mbox{a.e. in } \omega_0.
\end{equation}
\end{thm}
\begin{proof} 
By the results expressed in equations (2.15), (2.16), (2.17), (4.4) 
and lemma 2 we 
have to compute only the weak limit of the sequence $\{h_1^k\}$, where
\[
h_1^k=
a(x)\abs{\varepsilon (u^k (x))}^2.
\]
The sequence $\{p^k\cc\varepsilon(u^k)\}$ as being
weakly convergent in $L^1(\Omega)$ has to be equiintegrable according to the
Dunford-Pettis criterion of weak compactness in $L^1(\Omega)$. Since
$p^k\cc\varepsilon(u^k)=m^k\abs{\varepsilon(u^k)}^2+
\bigl(\mathcal{A}^++\psi^k \mathcal{A}^-)\cc\varepsilon(u^k)$, it can be easy to deduce that
$\{m^k\abs{\varepsilon(u^k)}^2\}$ is equiintegrable as well 
{(and so is $\{h_1^k\}$)}.
Thus one can
suppose that it is weakly
convergent in $L^1(\Omega)$, by passing to a subsequence, if
necessary,
so by Theorem 6.2, p. 97, \cite{PEDREGAL} we
see, that its weak limit is 
$\int\limits_{\mathbb{R}^{n\times n}}a\abs{\lambda }^2 d\nu_x (\lambda).$
\\
Now, from  the inequality $(R)$
\[
h\bigl({{x,}}\,
\varepsilon(u^k {{(x)}})\bigr)\leq
\frac{1}{2}m^k\abs{\varepsilon (u^k)}^2+
(\mathcal{A}^++\psi^k \mathcal{A}^-)\cc\varepsilon (u^k)
+\frac{1}{2}\mathcal{B}^k
\]
we are allowed to conclude that the sequence
$\{h\bigl({{x},}\,\varepsilon(u^k {{(x)}})\bigr)\}$
has the same property.
As shown in (Theorem 6.2, p. 97, \cite{PEDREGAL}) the weak limit is then a
function as just given on the right hand side of \eqref{P-7}.
\par
To show \eqref{P-71} it is enough to establish the strong convergence of
$\{\varepsilon(u^k)\}$ in $L^2(\omega_0;\allowbreak\RRR)$ (cf. Proposition 6.12,
p. 111, \cite{PEDREGAL}).
The elements of the sequence $\{\psi^k\}=\{\chi_b^k
-\chi_a^k\}$ take values $+1$ or $-1$. Thus the upper Kuratowski limit of the
sequence of singletons $\{\psi^k (x)\}$ (i.e. the set of limit points
of this sequence) is the set $\{-1,\,1\}$. By the Balder theorem (see 
\cite{VALADIER}) we see that $\psi^k\rightarrow 1$ strongly in
$L^1 (\omega_0^+ )$ and $\psi^k\rightarrow -1$ strongly in 
$L^1 (\omega_0^- )$
and we can suppose that $\psi^k\to 1$ 
a.e. in $\omega_0^{{+}}$ 
{{($\psi^k\to -1$ 
a.e. in $\omega_0^-$)}}
by passing to a subsequence, if necessary. Further, the
equiintegrability of $\{m^k\abs{\varepsilon(u^k)}^2\}$ implies that
$\{\abs{\varepsilon(u^k)}^2\}$ is also equiintegrable. By \lemref{lem-1} we have
\begin{equation*}
\int\limits_{\omega_0^{{+}}}p^k\cc\varepsilon(u^k)\,dx\to
\int\limits_{\omega_0^{{+}}}p\cc\varepsilon(u)\,dx=
\int\limits_{\omega_0^{{+}}}b\abs{\varepsilon(u)}^2\,dx+
\int\limits_{\omega_0^{{+}}}(\mathcal{A}^++\mathcal{A}^-)\cc\varepsilon(u)\,dx.
\end{equation*}
On the other hand,
\begin{equation*}
\int\limits_{\omega_0^{{+}}}p^k\cc\varepsilon(u^k)\,dx=
\int\limits_{\omega_0^{{+}}}b\abs{\varepsilon(u^k)}^2\,dx+
\int\limits_{\omega_{0k}^-}(a-b)\abs{\varepsilon(u^k)}^2\,dx
+\int\limits_{\omega_{0}^{{+}}}(\mathcal{A}^++\psi^k\mathcal{A}^-)\cc\varepsilon(u^k)\,dx,
\end{equation*}
where $\omega_{0k}^-=\{x\in\omega_0^+\colon
\psi^k(x)=-1\}$.
Thus taking into account that
\[
\int\limits_{\omega_{0}^{{+}}}(\mathcal{A}^++\psi^k\mathcal{A}^-)\cc\varepsilon(u^k)\,dx
\to\int\limits_{\omega_{0}^{{+}}}(\mathcal{A}^++\mathcal{A}^-)\cc\varepsilon(u)\,dx
\]
and
\[
\int\limits_{\omega_{0k}^-}(a-b)\abs{\varepsilon(u^k)}^2\,dx\to 0
\]
being a consequence of the equiintegrability of $\{\abs{\varepsilon(u^k)}^2\}$
and $\abs{\omega_{0k}^-}\to 0$,
we are led to
\[
\int\limits_{\omega_0^{{+}}}b\abs{\varepsilon(u^k)}^2\,dx\to
\int\limits_{\omega_0^{{+}}}b\abs{\varepsilon(u)}^2\,dx.
\]
Since, simultaneously, $\varepsilon(u^k)\rightharpoonup
\varepsilon(u)$ in $L^2(\omega_0^+ ;\mathbb{R}^{n\times n}_{sym})$, 
the desired strong convergence
results. Analogous reasoning holds for $\omega_0^-$.
The proof is complete.
\end{proof}
\begin{examp}
{{
Let $\omega_0=\{x\in \Omega\colon \chi_a(x)\chi_b(x)=0\}$. Without loss of generality one can suppose that
$\psi=1$ a.e. in $\omega_0$. Let $\omega_{0k}^-=\{x\in\omega_0\colon
\psi^k(x)=-1\}$. Since $\psi^k\to 1$ weak$^{\,\star}$ in
$L^\infty(\omega_0)$, we have
\begin{equation}\label{P-8}
2\abs{\omega_{0k}^-}=\int\limits_{\omega_0}(1-\psi_k)\,dx\to 0.
\end{equation}
Thus $\psi^k\to 1$ strongly in $L^1(\omega_0)$ (in fact, in
$L^p(\omega_0)$ for any $p\geq 1$). By the above theorem this means
that $\nu_x=\delta_{\varepsilon(u(x))}\quad\mbox{a.e. in } \omega_0.$
}}
\end{examp}
\begin{rem}
From \eqref{eq-1.28} and \eqref{P-7} it follows that
\begin{eqnarray*}
\lim\limits_{k\to\infty}\int\limits_{\Omega_0}a\abs{\varepsilon(u^k)}^2dx=
\int\limits_{\Omega_0}\int\limits_{\Rn}a\,\abs{\lambda}^2\,
d\nu_x(\lambda)\,dx=\int\limits_{\Omega_0}a\abs{\varepsilon(u)}^2dx
\\
+\theta\int\limits_{\Omega_0}\chi_a\chi_b\,a\abs{C-D}^2dx,
\end{eqnarray*}
giving rise to the formula that allows to calculate $\theta\in[0,1]$. Namely,
if we let
\begin{equation}\label{eq-d}
d:=\lim_{k\to\infty}\int\limits_{\Omega_0}a\abs{\varepsilon(u^k)}^2dx-
\int\limits_{\Omega_0}a\abs{\varepsilon(u)}^2dx,
\end{equation}
then from the equation
\begin{equation*}
d={\frac{\theta}{2}}\int\limits_{\Omega_0}\chi_a\chi_b\,a\abs{C-D}^2dx
\end{equation*}
we obtain
\begin{equation}\label{eq-dd}
\theta=\begin{cases} \displaystyle
\frac{2d}{\int\limits_{\Omega_0}\chi_a\chi_b\,a\abs{C-D}^2dx}&\mbox{if}\;
\int\limits_{\Omega_0}\chi_a\chi_b\,a\abs{C-D}^2dx>0\\[7mm]
0&\mbox{otherwise},
\end{cases}
\end{equation}
or equivalently
\begin{equation}\label{eq-ddd}
\theta=\begin{cases} \displaystyle
\frac{2\int\limits_{\Omega_0}\int\limits_{\Rn}a\,\abs{\lambda}^2\,
d\nu_x(\lambda)\,dx-\int\limits_{\Omega_0}a\abs{\varepsilon(u)}^2dx}
{\int\limits_{\Omega_0}\chi_a\chi_b\,a\abs{C-D}^2dx}&\mbox{if}\;
\int\limits_{\Omega_0}\chi_a\chi_b\,a\abs{C-D}^2dx>0\\[7mm]
0&\mbox{otherwise},
\end{cases}
\end{equation}
\end{rem}
\begin{rem}
It is worth to
point out that
the formulas
\eqref{eq-1.28}, \eqref{eq-1.30}, \eqref{eq-1.31} make possible
to express
the infimum of $(P)$  via \eqref{eq-dd} in terms of the limits $u$, $p$, $\chi_a$,
$\chi_b$, $d$ only.
On the other hand, the formula \eqref{P-7} expresses it in terms of the
parametrized Young measures $\{\nu_x(\cdot)\}$ which, in practice,
are  much more difficult to derive.
\end{rem}
\par
\vspace{1pc}
\noindent
\textbf{Acknowledgement.} I wish to express my gratitude to Professor
Krzysztof Che{\l}mi\'nski for his significant remarks and valuable suggestions

\bibliography{ref-main,ref-nan}

\providecommand{\bysame}{\leavevmode\hbox to3em{\hrulefill}\thinspace}

\end{document}